\documentclass[12pt]{article}

\usepackage{amsmath,amssymb,amsthm}
\chardef\bslash=`\\
\newtheorem{Thm}{Theorem}[section]
\newtheorem{Cor}[Thm]{Corollary}
\newtheorem{Lem}[Thm]{Lemma}

\newtheorem{Rque}[Thm]{Remark}
\newtheorem{Conj}{Conjecture}

\newtheorem{Exam}[Thm]{Example}

\title{Birational geometry in codimension 2 of symplectic resolutions}
\author{Baohua Fu}

\def\cit{{\mathbb C}}

\def\pit{{\mathbb P}}
\def\zit{{\mathbb Z}}
\def\0{{\mathcal O}}
\def\g{{\mathfrak g}}

\begin{document}
\maketitle
\begin{abstract}
We prove the conjecture that two projective symplectic resolutions for a symplectic variety $W$ are
related by Mukai's elementary transformations over $W$ in codimension 2 in the following cases:
(i). nilpotent orbit closures in a classical simple complex Lie algebra; (ii). some quotient
symplectic varieties.
\end{abstract}
\section{Introduction}

A symplectic variety is a complex
algebraic variety $W$, smooth in codimension 1, such that there exists a
regular symplectic form on its smooth part which can be extended  to
a global regular form on any resolution (see \cite{Bea}).  A resolution 
 $\pi: X \to W$ is called {\em symplectic} if 
 the lifted regular form on $X$ is non-degenerated everywhere. One can show that a resolution is symplectic if
and only if it is crepant.

Let $W$ be a symplectic variety and $\pi: X \to W$ a symplectic 
resolution. Assume that  
$W$ contains a smooth subvariety $Y$ such that the restriction of $\pi$  to $P: = \pi^{-1}(Y)$ 
makes $P$ a $\pit^n$-bundle over $Y$. 
If $codim(P) = n$, then we can blow up $X$ along $P$ and then blow down 
along the other direction, which gives another (proper) symplectic
resolution $\pi^+: X^+ \to W$, provided that $X^+$ remains in our 
category of algebraic varieties. The diagram $X \to W \leftarrow X^+$
is called  {\em  Mukai's elementary transformation} (MET for short) 
over $W$ with center $Y$. A {\em MET in codimension 2} is a diagram
which becomes a MET after removing subvarieties of codimension greater 
than 2. \vspace{0.4cm}

\begin{Conj}
Let $W$ be a symplectic variety which admits two projective symplectic resolutions 
$\pi: X \to W$ and $\pi^+: X^+ \to W$. Then the birational map
$\phi = (\pi^+)^{-1} \circ \pi: X --\to X^+$ is related by a sequence 
of METs over $W$ in codimension 2. 
\end{Conj}

Notice that since the two resolutions $\pi, \pi^+$ are both crepant, 
the birational map $\phi$ is isomorphic in codimension 1.

In \cite{HY} (Conjecture 7.3), Hu and Yau made an analogue conjecture 
for birational maps between smooth 
projective holomorphic symplectic varieties. Regarding the many alike
properties shared by hyperk\"ahler manifolds and projective symplectic resolutions,
we believe that their conjecture holds for projective symplectic 
resolutions.

This conjecture is true for four-dimensional symplectic varieties by the work of 
Wierzba and  Wi\'sniewski (\cite{WW}). The purpose of this note is to
provide more evidence for this conjecture. Using  results of Namikawa 
(\cite{Nam}), we prove the following:
\begin{Thm}
Let $\overline{\0}$ be a nilpotent orbit closure in a simple classical 
Lie algebra. Then any two (proper) symplectic resolutions for 
$\overline{\0}$ are connected  by a sequence of METs over 
$\overline{\0}$  in codimension 2. 
\end{Thm}

Then we consider the situation of quotient symplectic varieties. Let 
$V$ be a vector space and $G$ a finite subgroup of $GL(V)$,
then $G$ acts naturally on $T^*V$ and the quotient  $(T^*V)/G$ is a symplectic 
variety. In this case, we have the following:
\begin{Thm}
Let $G$ be a finite subgroup of $GL(V)$ such that for any codimension 2 
subspace $H \subset V$, the set $\{g \in G| V^g = H\}$
forms a single conjugacy class. Then for any two projective symplectic 
resolutions $\pi_i: Z_i \to (T^*V)/G, i=1, 2,$
the induced birational map $\phi: Z_1 --\to Z_2$ is isomorphic in 
codimension 2.
\end{Thm}

Note that there are only finitely many codimension 2 subspaces $H$ such that
the set $\{g \in G| V^g = H\}$ is non-empty. 
This theorem has an interesting corollary.

\begin{Cor}
Let $G$ be a finite subgroup of $GL(2)$  such that the elements $g \in 
G$ whose  eigenvalues are all different to $1$ form a single conjugacy 
class.
Then $ (T^* \cit^2)/G$ admits at most one projective symplectic 
resolution, up to isomorphisms.
\end{Cor}

In the proof, we obtain the following theorem (valid for
a general $G$), which could be helpful in further studies. 
 This also  gives an illustration of philosophy of the  McKay
correspondence: how the geometry of a symplectic resolution of $ (T^*V)/G$ is controlled by the
group $G$.
\begin{Thm}
Let $V$ be a vector space and $G < GL(V)$ a finite sub-group.
Suppose we have a projective symplectic resolution $\pi: Z \to (T^*V)/G$.
Then:

(i)(\cite{Ka1}). $V/G$ is smooth;

(ii).  $Z$ contains a Zariski open set $U$ which is isomorphic to $T^*(V/G)$;

(iii). the morphism $\pi: T^*(V/G) \to (T^*V)/G$ is the natural one (c.f.  section 3), 
which is independent of the resolution.
\end{Thm}

\section{Nilpotent orbits}
\subsection{Stratified Mukai flops}

Consider the nilpotent orbit $\0= \0_{[2^k, 1^{n-2k}]}$ in 
$\mathfrak{sl}_n$, where $2k \leq n$. The closure  $\overline{\0}$
admits exactly two symplectic resolutions given by $$T^*G(k,n) 
\xrightarrow{\pi}  \overline{\0} \xleftarrow{\pi^+} T^*G(n-k, n),$$
where $G(k,n)$ (resp. $G(n-k, n)$) is the Grassmannian of $k$ (resp 
$n-k$) dimensional subspaces
 in $\cit^n$. Let $\phi$ be the induced birational map 
$T^*G(k,n) --\to T^*G(n-k,n)$.

It is shown by Namikawa (\cite{Nam} Lemma 3.1) that when $2k < n$,  
$\pi$ and $\pi^+$ are both small  and the diagram is a flop.
This is the {\em stratified Mukai flop of type $A_{k,n}$.}
When $2k = n$, the birational map $\phi$ is an isomorphism. 
\begin{Lem} \label{typeA}
If $n \neq 2k+1$, then $\phi$ is an isomorphism  in codimension 2. If 
$n = 2k +1$, then $\phi$
is a MET over $\overline{\0}$ in codimension 2. 
\end{Lem}
\begin{proof}
The closure $\overline{\0}$ consists of orbits $\{\0_{[2^i, 1^{n-2i}]}  
\}_{0 \leq i \leq k}$.
The fiber of $\pi$ (resp. $\pi^+$) over a point in $\0_{[2^i, 
1^{n-2i}]} $ is isomorphic to $G(k-i, n-2i)$
(resp. $G(n-k-i, n-2i)$). By a simple dimension count, one shows that 
the complement of $\pi^{-1}(\0)$ (resp. $(\pi^+)^{-1}(\0)$ ) 
is of codimension greater than 2 when $n \neq  2k+1$, which proves that 
$\phi$ is isomorphic in codimension 2.

Now suppose that $n = 2k+1$. Let $Y$ be the nilpotent orbit 
$\0_{[2^{k-1}, 1^3]}$ and $P$ (resp. $P^+$) the preimage of $Y$ under $\pi$ (resp. 
$\pi^+$).
Then $P$ is  the subvariety $$\{([F], x) \in G(k, 2k+1) \times Y | 
Img(x) \subset F \subset Ker(x) \}$$ in
 $T^*G(k, 2k+1) \subset G(k, 2k+1) \times \overline{\0}.$  The induced 
map $P \to Y$ makes $P$ a $\pit^2$-bundle over $Y$.
Similarly $P^+$ is the subvariety $$\{([F^+], x) \in G(k+1, 2k+1) 
\times Y | Img(x) \subset F^+ \subset Ker(x) \}$$ in $T^*G(k+1, 2k+1)$. 
The map $P^+ \to Y$ makes $P^+$ a $\pit^2$-bundle over $Y$. 

Let $U = \0 \cup Y$, which is open in $\overline{\0}$. The  
complement of $\pi^{-1}(U)$ (resp.  $(\pi^+)^{-1}(U)$) is of codimension 
greater than
2. Notice that the $\pit^2$-bundle $P$ over $Y$ is the dual of the 
$\pit^2$-bundle $P^+$ over $Y$. One deduces that
the diagram $\pi^{-1}(U) \to U \leftarrow (\pi^+)^{-1}(U)$ is a MET over $U$ 
with center $P$, which concludes the proof.
\end{proof}

Notice that the precedent proof gives an explicit description of the center 
of the MET, which will be used later. 

Now we introduce the stratified Mukai flops of type $D$.  Let $\0$ be 
the orbit $\0_{[2^{k-1}, 1^2]}$ in $\mathfrak{so}_{2k}$, where
$k \geq 3$ is an odd integer. Let $G_{iso}^+(k), G_{iso}^-(k)$ be the 
two connected components of the orthogonal Grassmannian 
of $k$-dimensional isotropic subspace in $\cit^{2k}$ (endowed with a 
fixed non-degenerate symmtric form). Then we have two symplectic 
resolutions $T^*G_{iso}^+(k) \to \overline{\0} \leftarrow 
T^*G_{iso}^-(k)$. It is shown in \cite{Nam} (Lemma 3.2) that this diagram is a flop
and the two resolutions are both small.

 Let $\phi$ be the induced birational map from $T^*G_{iso}^+(k)$ to 
$T^*G_{iso}^-(k)$. Then a simple dimension count shows that:
\begin{Lem}\label{typeD}
$\phi$ is an isomorphism in codimension 2.
\end{Lem}
\subsection{$\g = \mathfrak{sl}_n$}
Let $\0$ be a nilpotent orbit in  $\mathfrak{sl}_n$ corresponding to 
the partition 
${\bf d} = [d_1, \cdots, d_k]$ and $x \in \0$. Let $(p_1, \cdots, p_s)$ 
be a sequence of integers such that
$d_i = \sharp \{j| p_j \geq i\}.$ Fix a flag $F:= \{ F_i\}$ of $\cit^n$ 
of type $(p_1, \cdots, p_s)$ such that
$x F_i \subset F_{i-1}$ for all $i$. Such a flag is called a 
{\em polarization} of $x$. Every nilpotent element has only finitely many
different polarizations.

Assume that $p_{j-1} < p_j$ for some $j$. Consider the map $\alpha: F_j 
\to F_j/F_{j-2}$.
The element $x$ induces $\bar{x} \in End(F_j/F_{j-2})$. 
We define a flag $F'$ by $F'_i = F_i$ if $i \neq j-1$ and  $F'_{j-1} = 
\alpha^{-1} (Ker(\bar{x}))$.
By Lemma 4.1 \cite{Nam}, $F'$ is again a polarization of $x$ with type 
$(p_1, \cdots, p_j, p_{j-1}, \cdots, p_s)$.

Let $P$ (resp. $P'$) be the stabilizer of $F$ (resp. $F'$) in $G= SL_n$. 
Then we obtain two symplectic resolutions
$T^*(G/P) \xrightarrow{\pi} \overline{\0} \xleftarrow{\pi'} T^*(G/P').$ 
Let $\phi: T^*(G/P) --\to T^*(G/P')$ be the induced birational map.
\begin{Lem}\label{sl}
(i) If $p_j \neq p_{j-1} +1$, then $\phi$ is isomorphic in codimension 
2;

(ii) If $p_j = p_{j-1} + 1$, then $\phi$ is a MET over $\overline{\0}$  
in codimneiosn 2.
\end{Lem}
\begin{proof}
To simplify the notations, set $r = p_{j-1}$ and $m = p_j + p_{j-1}$. 
Let $\bar{F}$ be the flag obtained from $F$
by deleting the subspace $F_{j-1}$, which is also obtained from $F'$ by 
the same manner. We denote by $\bar{P} \subset G$ the stabilizer
of $\bar{F}$.
Let $X$ be the subvariety in $G/\bar{P} \times \overline{\0}$ consiting 
of the points $(\bar{E}, y)$ such that (i) $y \bar{E}_i \subset 
\bar{E}_{i-1}$
for $i \neq j-1$ and $y \bar{E}_{j-1} \subset \bar{E}_{j-1}$; (ii) the 
induced map $\bar{y} \in End(\bar{E}_{j-1}/\bar{E}_{j-2})$
satisfies $\bar{y}^2 =0$ and $rank(\bar{y}) \leq r$.

The projection to the second factor of $G/\bar{P} \times \overline{\0}$  induces a morphism $pr: X \to 
\overline{\0}$. By the proof of Lemma 4.3 \cite{Nam},
the resolutions $\pi, \pi'$ factorize through the map $pr$, which gives 
a diagram $T^*(G/P) \xrightarrow{\mu} X \xleftarrow{\mu'} T^*(G/P'). $
By Lemma 4.3 \cite{Nam}, this diagram is locally a trivial family of 
stratified Mukai flops of type $A_{r,m}$. By Lemma \ref{typeA}, 
if $m \neq 2r + 1$, then $\phi$ is isomorphic in codimension 2, which 
proves claim (i).

Now assume that $m = 2r+1$. Let $Y$ be the subvariety in $X$ consists 
of the points $(\bar{E}, y)$ such that the induced map
$\bar{y} \in End(\bar{E}_{j-1}/\bar{E}_{j-2})$ has rank $r-1$. By the 
proof of Lemma \ref{typeA} and Lemma 4.3 \cite{Nam},
the diagram  $T^*(G/P) \xrightarrow{\mu} X \xleftarrow{\mu'} T^*(G/P')$ 
is a MET over $X$ in codimension 2 with center $Y$.

Let ${\bf d'}$ be the partition of $n$ given by (possibly we need to re-order these parts):  
$$
d'_i =  \begin{cases}  d_i,  \ & \text{if} \ i  \neq r,  r+2 \\  d_r -1,  \ & \text{if} \ i =r \\
  d_{r+2} + 1, \ & \text{if} \ i=r+2. \end{cases}
$$

Then one can verify that the morphism
 $pr: X \to \overline{\0}$ maps $Y$ isomorphically to the nilpotent 
orbit
$\0_{{\bf d'}}$, which shows that the diagram $T^*(G/P) 
\xrightarrow{\pi} \overline{\0} \xleftarrow{\pi'} T^*(G/P')$ is a MET in codimension 2
over $\overline{\0}$ with center $\0_{{\bf d'}}$. 
\end{proof}

Notice that the precedent  proof  gives an explicit way to find out the 
MET center in $\overline{\0}$. Here we give an example.
\begin{Exam}
(Example 4.6 \cite{Nam}).
Let $\0 = \0_{[3,2,1]} \subset \mathfrak{sl}_6$ and $x \in \0$. Then 
$x$ has six polarizations $P_{\sigma(1), \sigma(2), \sigma(3)}$
of flag type $( \sigma(1), \sigma(2), \sigma(3))$, where $\sigma$ are 
permutations. Let $Y_{i,j,k} = T^*(SL_6/P_{i,j,k})$,
which gives a symplectic resolution for $\overline{\0}$. Then $Y_{321} 
--\to Y_{231}$ is a MET in codimension 2 with center $\0_{[3, 1^3]}$;
$Y_{231} --\to Y_{213}$ is isomorphic in codimension 2; $Y_{213} --\to 
Y_{123}$ is a MET in codimension 2 with center $\0_{[2^3]}$ and so on.
If a center appears twice in a sequence, then it is not really 
a MET center. For example, the birational map
$Y_{321} --\to Y_{132}$ is  a MET in codimension 2 with center 
$\0_{[2^3]}$, but over the orbit $\0_{[3, 1^3]}$, it is an isomorphism.
\end{Exam}

\subsection{$\g = \mathfrak{so}(V)$ or $\mathfrak{sp}(V)$}

Let $V$ be an $n$-dimensional vector space endowed with a 
non-degenerate bilinear symmetric (resp. anti-symmetric) form for  $\g = 
\mathfrak{so}(V)$
(resp.  $\g = \mathfrak{sp}(V)$).
Let $\epsilon =0$ if $\g = \mathfrak{so}(V)$ and $\epsilon = 1$ if $\g 
= \mathfrak{sp}(V)$.

Let $P_\epsilon(n)$ be the set of partitions ${\bf d}$ of $n$ such that 
$\sharp \{i| d_i=l\}$ is even for every integer $l$ with 
$l \equiv \epsilon$ (mod 2). These are exactly those partitions which 
appear as the Jordan types of nilpotent elements of 
$\mathfrak{so}(V)$ or of $\mathfrak{sp}(V)$.  Let $q$ be a non-negative 
integer such that $q \neq 2$ if $\epsilon =0$.
Define $Pai(n,q)$ to be the set of partitions ${\bf e}$ of  $n$ such 
that $e_i \equiv 1$ (mod 2) if $i \leq q$ and 
$e_i \equiv 0$(mod 2) if $i > q$. For ${\bf e} \in Pai(n, q)$, let 
$$I({\bf e}) = \{j | j \equiv n+1 \text{(mod 2)}, e_j \equiv \epsilon \text{(mod 2)}, 
e_j \geq e_{j+1}+2 \}.$$

The Spaltenstein map $S: Pai(n,q) \to P_\epsilon(n)$ is defined as 
$$S({\bf e}) _j = \begin{cases} e_j - 1,  &\text{if} \ j \in I({\bf e}) \\
 e_j + 1, &\text{if} \ j-1 \in I({\bf e}) \\  e_j, &\text{otherwise}. \end{cases} $$

It is proved in \cite{Hes}
that for a nilpotent element of type ${\bf d}$, its polarization types 
are determined by $S^{-1}({\bf d})$. 
For a sequence of integers $(p_1, \cdots, p_k)$, we define ${\bf e}= 
ord(p_1, \cdots, p_k)$ to be the partition given by
$e_i = \sharp \{j|p_j \geq i\}$.

Let $\0$ be a nilpotent orbit of type {\bf d} in $\g$ and $x\in \0$. 
Let $(p_1, \cdots, p_k, q, p_k, \cdots, p_1)$ be a sequence of
integers such that ${\bf e}: = ord(p_1, \cdots, p_k, q, p_k, \cdots, 
p_1)$ is in $Pai(n, q)$ and $S({\bf e}) = {\bf d}.$
Let $F$ be an isotropic flag (i.e. $F_i^\bot = F_{2k+1-i}, \forall i$)  
in $V$ of type  $(p_1, \cdots, p_k, q, p_k, \cdots, p_1)$
such that $x F_i \subset F_{i-1}$.

Assume that $p_{j-1} < p_j$ for some $j$.  Consider the map $\alpha: 
F_j \to F_j/F_{j-2}$.
The element $x$ induces $\bar{x} \in End(F_j/F_{j-2})$. 
We define another flag $F'$ by $F'_i = F_i$ if $i \neq j-1, 2k+2-j$, 
$F'_{j-1} = \alpha^{-1} (Ker(\bar{x}))$ and
$F'_{2k+2-j} = (F'_{j-1})^\bot.$ By Lemma 4.2 \cite{Nam}, $F'$ is again 
a polarization of $x$.
We denote by $P$ (resp. $P'$) the stabilizer of $F$ (resp. $F'$). Then 
we obtain two symplectic resolutions
$T^*(G/P) \xrightarrow{\pi} \overline{\0} \xleftarrow{\pi'} T^*(G/P').$ 
Let $\phi$ be the induced birational map
from $T^*(G/P) $ to $T^*(G/P')$. 
\begin{Lem}\label{sop}
(i). If $p_j \neq p_{j-1} + 1$, then $\phi$ is isomorphic in 
codimension 2;

(ii). If $p_j = p_{j-1} + 1$, then $\phi$ is a MET in codimension 2 
over $\overline{\0}$.
\end{Lem}

The proof goes along the same line as that in Lemma \ref{sl}. The 
difference is the definition of the partition ${\bf d'}$ in the proof of 
(ii).
Here we have $r = p_{j-1}$ and $p_j = r+1$. Let ${\bf e'}$ be the 
partition (after re-ordering if necessary)  defined by 

$$e'_j = \begin{cases} e_j, & \text{if} \ j \neq r, r+2 \\ 
 e_r -2, & \text{if} \ j=r \\ e_{r+2} + 2, & \text{if} j = r+2. \end{cases} $$

Then ${\bf e'} \in Pai(n,q)$. 
Now  we should define ${\bf d'}  = S({\bf e'})$. In this case, $\phi$ 
is a MET in codimension
2 over $\overline{\0}$ with center $\0_{{\bf d'}}$.
  
\begin{Exam} (Example 4.7 \cite{Nam}).
 Let $\0 = \0_{[4^2, 1^2]}$ be the nilpotent orbit in 
$\mathfrak{so}_{10}$. Take an element $x \in \0$, then $x$ has 
four polarizations $P_{3223}^+, P_{3223}^-, P_{2332}^+, P_{2332}^-$. 
Let $Y_{3223}^+ = T^*(G/P_{3223}^+)$ and so on. Then 
$Y_{3223}^+  --\to Y_{2332}^+ $ is a MET in codimension 2 over 
$\overline{\0}$ with center $\0_{[3^2,2^2]}$.
\end{Exam}

\subsection{Proof of theorem 1.1}

Let $\0$ be a nilpotent orbit in a classical simple Lie algebra $\g$. 
By \cite{Fu}, every (proper) symplectic resolution for
$\overline{\0}$ is of the form $T^*(G/P) \to \overline{\0}$ for some 
polarization $P$ of $\0$. 
Assume  that we have two symplectic resolutions $T^*(G/P_i) \to 
\overline{\0}, i=1,2$, then by the proof of Theorem 4.4 \cite{Nam},
we can reach $T^*(G/P_2) \to \overline{\0}$ from $T^*(G/P_1) \to 
\overline{\0}$ by using the operations in section 2.2 and 2.3,
possibly by using another operation which is a locally trivial family of 
stratified Mukai flops of type $D$ (thus isomorphic in codimension 2
by Lemma \ref{typeD}). Now Lemma \ref{sl} and \ref{sop} give the 
theorem. 

\section{Quotient singularities}
Let $V$ be an $n$-dimensional vector space and $G$ a finite subgroup of 
$GL(V)$. Then $G$ acts naturally on $T^*V \simeq V \oplus V^*$, 
preserving
the symplectic form on $T^*V$. By \cite{Bea}, $W: = (T^*V)/G$ is a 
symplectic variety. We will study projective symplectic resolutions of $W$.

Consider the $\cit^*$-action on $T^*V$ defined by $\lambda (v, v') = 
(v, \lambda v')$ for $\lambda \in \cit^*$ and $(v, v') \in V \oplus V^*$.
This action commutes with the action of $G$, thus we obtain an action 
of $\cit^*$ on $W$. The fixed point set under this action is identified 
with 
$V/G$. If we denote by $\omega$ the symplectic form on the smooth part 
of $W$, then $\lambda^* \omega = \lambda \omega$ for all $\lambda \in \cit$.

Let $Z \xrightarrow{\pi} W$ be a projective symplectic resolution. By \cite{Ka1}, 
the $\cit^*$-action on $W$ lifts to a $\cit^*$-action on $Z$ in such
a way that $\pi$ is $\cit^*$-equivariant. If we denote by $\Omega$ the 
symplectic form extending $\pi^* \omega$ to the whole of $Z$, 
then $\lambda^* \Omega = \lambda \Omega$  for any $\lambda \in \cit^*$. 
Let $Z^{\cit^*}$ be the  points of $Z$ fixed by the $\cit^*$-action.
\begin{Lem} 
There exsits a connected component $Y$ of  $Z^{\cit^*}$ such that $\pi: 
Y \to V/G$ is an isomorphism. In particular, $V/G$ is smooth.
\end{Lem}

This is proved by  Kaledin  (the proof of Theorem 1.7 \cite{Ka1}). 
The following is proved in \cite{Fu2}, but since
\cite{Fu2} will never be published, we include the proof here.
Let $U = \{z \in Z| \lim_{\lambda \to 0} \lambda \cdot z \in Y \}.$
\begin{Lem}\label{open}
$U$ is isomorphic to $T^*Y$, and the induced $\cit^*$-action on $T^*Y$ 
is the natural action: $\lambda (y, v') = (y, \lambda v')$,
for any $y \in Y$, $v' \in T^*_yY$ and $\lambda \in \cit^*$.
\end{Lem}
\begin{proof}
 The precedent lemma shows that $Y$  is isomorphic to $V/G$, in 
particular $dim(Y) = n$.
For any point $y \in Y$, the action of $\cit^*$ on $Z$ induces a 
weight decomposition 
$$T_yZ = \oplus_{p \in \zit} T_y^p Z, $$ where $T_y^p Z = \{ v \in T_y 
Z |  \lambda_* v = \lambda^p v \}$, and $T_yY$ is identified with 
$T_y^0Z$.
The relation $\lambda^* \Omega = \lambda \Omega$ gives a duality 
between $T_y^p(Z)$
 and $T_y^{1-p}(Z)$. In particular,
$dim(T_y^1Z) = dim(T_y^0Z) = dim Y = n$, so $T_y^pZ =0$ for all $p \neq 
0,1$, which gives a decomposition $T_yZ = T_yY \oplus T_y^1Z$. 
Furthermore $Y$ is 
Lagrangian with respect to $\Omega$.

 By the work of Bialynicki-Birula (\cite{BB}), the decomposition $T_yZ 
= T_yY \oplus T_y^1Z$ shows that $U$ is a vector 
bundle of rank $n$ over $Y$, so $U$ is identified with the total space 
of the normal bundle $N$ of $Y$ in $Z$. 
Now we establish an isomorphism between $N$ and $T^*Y$ as follows. 
Denote by $\Omega_{can}$ the canonical symplectic structure on $T^*Y$.
Take a point $y \in Y$, and a vector $v \in N_y$.
Since $Y$ is Lagrangian in the both symplectic spaces, there exists a 
unique vector $w \in T_y^*Y$ such that
 $\Omega_y(v,u) = \Omega_{can,y}(w,u)$ for all $u \in T_yY$.
We define the map $i: N \rightarrow T^*Y$ to be $i(v)=w$. It is clear 
that $i$ is a $\cit^*$-equivariant isomorphism.
\end{proof}

Now we will study in more detail the morphism: $\pi: T^*Y \simeq 
T^*(V/G) \to (T^*V)/G.$ We denote by $p: V \to V/G$ the natural 
projection and $p_*: TV \to T(V/G)$ the induced tangent morphism.

We define a morphism  $\tilde{p}: T^*(V/G) \to (T^*V)/G $ as follows:
take a point $[x] \in V/G$ and a co-vector $\alpha \in T_{[x]}^*(V/G)$. We define a 
co-vector $\beta \in T_x^*V$ by
$<\beta, v> = <\alpha, p_*(v) > $ for all $v \in T_x V$. Then we put
$\tilde{p} ([x], \alpha) = [x, \beta]$. 
\begin{Lem}
$\tilde{p}$ is well-defined.
\end{Lem}
\begin{proof}
Let $y = gx$ with $g \in G$. We consider $\beta' \in T_y^*V$ defined by 
$<\beta', w> = <\alpha, p_*(w)>$
for all $w \in T_yV = g_* T_xV$. Then $w = g_* v$ for some $v \in T_x 
V$.  Now we have 
$$<\beta', g_*v> = <\alpha, p_*g_*(v)> =  <\alpha, p_*(v)>  = <\beta, v>$$
for all $v \in T_xV$, which gives  $\beta' = g^*\beta$. Then $[y, 
\beta'] = [gx, g^*\beta] = [x, \beta] $ in $(T^*V)/G$.
\end{proof}
Notice  that $p: V \to V/G$ is \'etale  at a point $x \in V$ if and only if 
 the stabilizer $G_x$ of $x$ in $G$
 is trivial, thus $\tilde{p}_{[x]}: T^*_{[x]}(V/G) \to (T^*_xV)/G$
is an isomorphism if and only if $G_x$ is trivial.
Furthermore $\tilde{p}$ induces an identity on the zero section $V/G$.

\begin{Lem}\label{same}
The morphism $\pi: T^*Y \simeq T^*(V/G) \to (T^*V)/G $ coincides with 
the morphism $\tilde{p}$.
\end{Lem}
\begin{proof}
Let $V_0 = V - \cup_{g \neq 1} V^g$, on which $G$ acts freely. 
$\tilde{p}$ induces an isomorphism 
$\tilde{p}_0: T^*(V_0/G) \to (T^*V_0)/G$.  We will show that $f:= 
\tilde{p}_0^{-1} \circ \pi|_{T^*(V_0/G)}:  T^*(V_0/G)  \to  T^*(V_0/G) $ is an
identity. Notice that $f$ is an identity over the zero section $V_0/G$. 
Furthermore $f$ is $\cit^*$-equivariant and it preserves the natural 
symplectic
form on $T^*(V_0/G)$.

Take a point $[x] \in V_0/G$ and $\alpha \in T_{[x]}^*(V_0/G)$. We 
consider $\alpha$ as a vector in $T_{[x]}(T^*(V_0/G)) $, then
$ f_* (\alpha) = \frac{d}{d \lambda}|_{\lambda =0} f(\lambda \alpha) = 
\frac{d}{d \lambda}|_{\lambda =0} \lambda f(\alpha) = f(\alpha).$
Now since $f$ preserves the symplectic form, we have $<\alpha, v> = 
<f_*(\alpha), f_*(v)> = <f(\alpha), v> $ for all $v \in T_{[x]}(V_0/G)$,
which gives $f(\alpha)= \alpha$, thus $f$ is an identity on $T^*(V_0/G)$.
Furthermore, the birational map  $\tilde{p}_0^{-1} \circ \pi$ is an identity
on $V/G$.

Now we need only to show that $\tilde{p}_0^{-1} \circ \pi$ is a morphism. Suppose it is not defined at 
some point $([x], \alpha)$ with $[x] \in V/G$ and $\alpha \in T_{[x]}^*(V/G)$.
Since $\pi$ and $\tilde{p}$ are both $\cit^*$-equivariant, $\tilde{p}_0^{-1} \circ \pi$ is not defined
at $\lambda \cdot ([x], \alpha) = ([x], \lambda \alpha)$ for all $\lambda \in \cit^*$.
Let $\lambda \to 0$, then one gets that $\tilde{p}_0^{-1} \circ \pi$ is not defined at the point 
$([x], 0) \in V/G$, a contradiction.
\end{proof}
\begin{Thm}
Let $G$ be a finite subgroup of $GL(V)$ such that for any codimension 2 
subspace $H \subset V$, the set $\{g \in G| V^g = H\}$
forms a single conjugacy class. Then for any two projective symplectic 
resolutions $\pi_i: Z_i \to (T^*V)/G, i=1, 2,$
the induced birational map $\phi: Z_1 --\to Z_2$ is isomorphic in 
codimension 2.
\end{Thm}
\begin{proof}
Let $g$ be an element such that $H:=V^g$ is of codimension 2. Let $W_g 
= (T^*H)/C(g)$, where $C(g)= Stab(H)/Cent(H)$ is 
the quotient of the subgroup $Stab(H)$   of elements $h \in G$ which 
preserve $H$ by the subgroup $Cent(H)$
of elements  $h \in Stab(H)$ which acts as identity on $H$.
$W_g$ is of codimension 4 in $W$, whose preimage by $\pi_i$ is of 
codimension 2, by the semi-smallness of projective 
symplectic resolutions (see \cite{Ka1}).

By the McKay correspondence (see \cite{Ka2}),  the number of 
codimension 2 components in $\pi_i^{-1}(W_g)$ equals to the 
the number of conjugacy classes of $\{h \in G| V^h = H\}$, which is 1 
by the hypothesis. 

Let $U_i$ be the Zariski open subset in $Z_i$ given by Lemma \ref{open} 
and $q_i: U_i \simeq T^*Y_i \to Y_i$ the natural projection.
Let $F_i \subset Y_i \simeq V/G$ be the subvariety $H/C(g)$. By Lemma 
\ref{same} and the explicite description of $\tilde{p}$, we have
 $q_i^{-1}(F_i) \subseteq \pi_i^{-1}(W_g)$.
Notice that $q_i^{-1}(F_i)$ is of codimension 2,
so its closure is the unique codimension 2 component of $\pi_i^{-1}(W_g)$,  then the complement of $U_i 
\cap  \pi_i^{-1}(W_g)$ is of codimension at least 3 in
$\pi_i^{-1}(W_g)$.   
By Lemma \ref{same}, the birational map $\phi$ induces an isomorphism  
between $U_1$ and $U_2$. 
The above arguments show that $\phi$ gives an isomorphism between the 
generic points of the unique codimension 2
components in $\pi_i^{-1}(W_g)$, for all $g \in G$ such that $V^g$ is 
of codimension 2 in $V$, thus $\phi$ is isomorphic in
codimension 2.
\end{proof}
\begin{Cor}
Let $G$ be a finite subgroup in $GL(2)$  such that the elements $g \in 
G$ whose  eigenvalues are all  different to $1$ form a single conjugacy 
class.
Then $ (T^* \cit^2)/G$ admits at most one projective symplectic 
resolution, up to isomorphisms.
\end{Cor}
\begin{proof}
By the hypothesis,
the set $\{g \in G| V^g = 0\}$ forms a single conjugacy class, thus
the two symplectic resolutions are isomorphic in codimension 2.
By the work of \cite{WW}, any two projective symplectic resolutions for 
$ (T^* \cit^2)/G$
are connected by Mukai flops, so there is no flop at all, thus  the two resolutions are isomorphic.
\end{proof}
\begin{Rque}
This theorem and the corollary can be regarded as an extension of Corollary 2.4 in \cite{FN}
(see also theorem 1.9 \cite{Ka1}), where it is proved that if the complex reflections in $G$ 
 form 
a single conjugacy class, then
the quotient admits at most one projective symplectic resolution (up to isomorphisms). 
\end{Rque}

\begin{Exam}
Here we give an example to show that our assumption on $G$ cannot be removed.
Let $x,y$ be the coordinates of $\cit^2$ and $G$ the subgroup of $GL(2)$
generated by the two elements: $\sigma(x,y) = (x, -y)$ and $\tau(x,y) = (y,x)$.
Then $G$ is the dihedral group of order 8, which acts on $T^*(\cit^2)$ (with
coordinates $x, z, y, w$) as follows:
 $$   \sigma(x, z, y, w) = (x, z, -y, -w), \quad \tau (x, z, y, w) = (y,  w, x,  z).$$

The quotient $T^*(\cit^2)/G$ is isomorphic to $Sym^2(\bar{S})$, where $\bar{S} = \cit^2/{\pm 1},$
which posseds exactly two non-isomorphic symplectic resolutions (for details, see Example 2.7 \cite{FN}),
one is connected to the other by a Mukai's elementary transformation. 
\end{Exam}
\begin{Exam} Here we give an example
to show that our corollary covers some situations where Corollary 2.4 in
\cite{FN} does not apply. 
Let $G$ be the subgroup of $GL(2)$ generated by the following two 
elements:
$$  \sigma: (x, y) \mapsto (-x, y); \quad \tau: (x,y) \mapsto (x, -y). 
$$
There is only one element $\sigma \circ \tau$ whose eigenvalues are 
different to $1$.
By the precedent Corollary, two projective symplectic resolutions for  $(T^* 
\cit^2)/G$ are isomorphic.

In fact, $(T^* \cit^2)/G$ is isomorphic to the product of two 
$A_1$-singlarities, thus it admits a
symplectic resolution $T^*\pit^1 \times T^*\pit^1 \to (T^* \cit^2)/G$. 
Notice that the 
unique 2-dimensional fibre is isomorphic to $\pit^1 \times \pit^1$, 
thus no MET over $(T^* \cit^2)/G$ can be performed. By \cite{WW},
this is the unique (up to isomorphisms) projective  symplectic resolution for 
$(T^* \cit^2)/G$.
\end{Exam}

\quad \\
Labortoire J. Leray, Facult\'e  des sciences \\
2, Rue de la Houssini\`ere,  BP 92208 \\
F-44322 Nantes Cedex 03 - France\\
\quad \\
baohua.fu@polytechnique.org
\end{document}